 \magnification=1200
\input amssym.def
\input amssym.tex
 \font\newrm =cmr10 at 24pt
\def\bul{\raise .9pt\hbox{\newrm .\kern-.105em } }

 \def\fr{\frak}
 \font\sevenrm=cmr7 at 7pt

 \baselineskip 20pt
 
 \def\h{\hbox{ }}

 \def\u{{\fr u}}

 \def\n{{\fr n}}
 \def\a{{\fr a}}

 \def\ss{{\fr s}}
 
 \def\b{{\fr b}}
 
 \def\hh{{\fr h}}

 \def\g{{\fr g}}

 \def\<{\le}
 \def\>{\ge}

 \def\s{{\h\subset\h}}
 
 \def\vs{\vskip }

 \def\mapright#1
  {\smash{\mathop
  {\longrightarrow}
  \limits^{#1}}}

 \def\kk#1{{\kern .4 em} #1}
 \def\vs{\vskip 1pc}

\hsize = 31pc
\vsize = 45pc
\overfullrule = 0pt
\rm

\centerline{\bf $\hbox{\rm Cent}\,U(\n)$ and a construction of
Lipsman--Wolf}
\vskip 1.5pc \centerline{Bertram Kostant}\vskip 1.5pc
{\bf Abstract:} Let
$G$ be a complex simply-connected semisimple Lie group and let $\g=
\hbox{\rm Lie}\,G$. Let $\g = \n_- +\hh + \n$ be a triangular
decomposition of $\g$. The authors in [LW] introduce a very nice
representation theory  idea for the construction of certain elements in
$\hbox{\rm cent}\,U(n)$. A key lemma in [LW] is incorrect but the idea 
is in fact valid. In our paper here we modify the construction so as to yield
the desired elements in $\hbox{\rm cent}\,U(\n)$. \vskip .5pc

{\bf Key words:} cascade of orthogonal roots, invariant theory
\vskip 4pt
{\bf MSC 2010 codes:} 22Exx 

\vskip 6pt

\centerline{\bf 1. Introduction}\vskip .5pc {\bf 1.1.} Let
$G$ be a complex simply-connected semisimple Lie group and let $\g=
\hbox{\rm Lie}\,G$. Let $\g = \n_- +\hh + \n$ be a triangular decomposition of
$\g$. If $\ss$ is a complex Lie algebra, then $U(\g)$ and $S(\g)$
will denote respectively the enveloping algebra and symmetric
algebra of $\g$. 

Some time ago I introduced what is presently
referred to as the cascade of orthogonal roots. Using the cascade, 
Anthony Joseph and I independently obtained, with very different methods,
a number of structure theorems of $S(\n)^{\n}$ (or equivalently
$\hbox{\rm cent} \,U(\n)$). The cascade is also used in [LW]. The present
paper deals with a neat, interesting representation-theoretic idea in [LW]
for constructing certain elements in $S(\n)^{\n}$. Basically the
construction begins with the linear functional $f$ on $U(\g)$ obtained,
as a matrix unit, involving the highest and lowest weight vectors of the
irreducible representation $\pi_{\nu}$ of $G$ with highest weight
$\nu$. Without unduly detracting from the idea we point out here that a
key lemma in [LW] (Lemma 3.7) is incorrect as it standards. A
counterexample is given in our present paper. Next we show that the
construction can be modified so as to produce a correct result. The
modification is of independent interest in that it introduces the notion
of what we call the codegree of a linear functiomal like $f$. If the
codegree of $f$ is
$k$,  then one obtains a harmonic element $f_{(k)}$ of degree $k$
in $S(\g)$. We then go on to show that $f_{(k)}$ is the desired element
in
$S(\n)^{\n}$. The main result is Theorem 2.8.

\vskip 1pc

\centerline{\bf 2. Lipsman--Wolf construction}\vskip .5pc
{\bf 2.1.} Let $\g$ be a complex semisimple Lie algebra and let
$G$ be a simply-connected complex Lie group such that $\g= \hbox{\rm Lie}\,G$. 
Let
$\ell = \hbox{\rm rank}\,\g$ and let ${\cal B}$ be the Killing form $(x,y)$ on
$\g$. If $\u$ is a complex vector space, $S(\u)$ will denote the
symmetric (graded) algebra over $\u$. ${\cal B}$ extends to a
nonsingular symmetric bilinear form, still denoted by ${\cal B}$, on
$S(\g)$ where, if
$x,y\in \g$ and $m,n\in {\Bbb N}$, then $(x^n,y^m) = 0$ if $m\neq n$
and  $$(x^n,y^n) = n!(x,y)^n.$$
 One may then identify $S(\g)$ with the algebra of
polynomial functions on $\g$ where $$x^n(y) = (x,y)^n.$$ Also
$u\mapsto \partial_u$ defines an isomorphism of $S(\g)$ with the
algebra of differential operators on $\g$ with constant
coefficients. If $u,y,w\in S(\g)$ one readily has $$(u,vw) =
(\partial_vu,w). \eqno (2.1)$$ The symmetric algebra $S(\g)$
becomes a degree-preserving $G$-module by extending, as a group of 
automorphisms, the adjoint action of
$G$ on $\g$. Let $J= S(\g)^G$. Then one knows (Chevalley) that
$J={\Bbb C}[p_1,\ldots,p_{\ell}]$ is a polynomial ring with
homogeneous generators $p_1,\ldots,p_{\ell}$. A polynomial $f\in
S(\g)$ is called harmonic if $\partial_{p_1}f =0, \,i=1,\dots,\ell.$
Let $H\s S(\g)$ be the graded subspace of all harmonic polynomials
in $S(\g)$. Then $H$ is a $G$-submodule of $S(\g)$ and one knows
$$S(\g)= J\otimes H\eqno (2.2)$$ is a $G$-module decomposition of
$S(\g)$. Explicitly, for any $k\in {\Bbb N}$, $H^k$ is given by
$$H^k= \hbox{The $\Bbb C$ span of}\,\,\{w^k\mid\,\,
\hbox{where $w\in \g$ is nilpotent.}\} \eqno (2.3)$$

\vs {\bf 2.2.}
Let $\hh$ be a Cartan subalgebra of $\g$ and let $\Delta$ be the
set of roots of $(\hh,\g)$. Also for any $\varphi\in \Delta$ let
$e_{\varphi}\in \g$ be a corresponding root vector. If $\ss\s\g$
is any subspace which is stable under $\hbox{\rm ad}\,\hh$ let $\Delta(\ss)= 
\{\varphi\in \Delta\mid\,e_{\varphi}\in \ss\}$. Let $\b$ be a Borel
subalgebra of $\g$  which contains $\hh$ and let $\n$ be the
nilradical of $\b$. A system of positive roots $\Delta_+$ in
$\Delta$ is chosen so that $\Delta_+ = \Delta(\n)$. Let $\b_-$ be
the Borel subalgebra containing $\hh$ which is ``opposite" to
$\b$. One then has a triangular decomposition
$$\g = \n_- + \hh +\n$$ where $\n_-$ is the nilradical of $\b_-$. If
$\Delta_- = \Delta(\n_-)$, then of course $\Delta_- = -\Delta_+$. 

Let $L$ be the lattice  of integral linear forms on
$\hh$ and let $L_o$ be the sublattice of integral linear forms which are
also in the root lattice. If $M$ is any locally finite $G$-module and
$\mu\in L$, then $M(\mu)$ will denote the $\mu$-weight space. Let
$\hbox{\rm Dom}(L)$ (resp.\ $\hbox{\rm Dom} (L_o))$ be the set of dominant
elements in $L$ (resp.\
$L_o$). For any $\lambda\in \hbox{\rm Dom}(L)$ let
$\pi_{\lambda}:G\to \hbox{\rm Aut}\,V_{\lambda}$ be some fixed irreducible
finite-dimensional representation with highest weight
$\lambda$. One knows that $$V_{\lambda}(0) \neq 0 \iff \lambda\in
\hbox{\rm Dom}(L_o).\eqno (2.4)$$ 
If $M$ is a locally finite $G$-module,  let 
$M_{(\lambda)}$ be the primary component of $M$ corresponding to $\lambda$. It
is obvious that if $\mu \in L$ and $S(\g)(\mu) \neq 0$, then $\mu\in
L_o$. In particular if $\lambda\in \hbox{\rm Dom}(L)$ and
$S(\g)_{(\lambda)}\neq 0$,  then
$\lambda \in  \hbox{\rm Dom}(L_o)$. On the other hand if $\lambda\in 
 \hbox{\rm Dom}(L_o)$, 
then (2.2) readily implies $$S(\g)_{(\lambda)} = J\otimes
H_{(\lambda)}.\eqno (2.5)$$ For any $\lambda\in \hbox{\rm Dom}(L_o)$ let
$$\ell(\lambda) = \hbox{\rm dim}\,V_{\lambda}(0).\eqno (2.6)$$ Then one knows
$$\eqalign{\hbox{dim Hom}_G(V_{\lambda}, H)& =\hbox{Hom}_G(V_{\lambda},
H_{(\lambda)})\cr & = \ell(\lambda).\cr}\eqno (2.7)$$ Let $\sigma_i,\,
i=1,\ldots,\ell(\lambda)$, be a basis of $\hbox{\rm Hom}_G(V_{\lambda},
H_{(\lambda)})$ and put
$H_{\lambda,i} = \sigma_i(V_{\lambda})$ so that one has a complete
reduction of $H_{(\lambda)}$,
$$ H_{(\lambda)}= \sum_{i=1}^{\ell(\lambda)} H_{\lambda,i}\eqno (2.8)$$
into a sum of irreducible components. Furthermore we can choose the
$\sigma_i$ so that
$H_{\lambda,i}$ is homogeneous for all $i$. In fact there is a unique
nondecreasing sequence of integers,
$m_i(\lambda),\,i=1,\ldots, m_{\ell(\lambda)}(\lambda)$, which are referred
to as  generalized exponents such that $$H_{\lambda,i}\s
H^{m_i(\lambda)}.\eqno (2.9)$$ Moreover the maximal generalized
exponent, $m_{\ell(\lambda)}(\lambda)$, occurs with multiplicity 1.  That
is $$m_i(\lambda)< m_{\ell(\lambda)}(\lambda)\eqno (2.10)$$ for any $i<
\ell(\lambda)$ and the maximal generalized
exponent is explicitly given by $$ m_{\ell(\lambda)}(\lambda) =
\sum_{i=1}^{\ell} k_i\eqno (2.11)$$ where $\lambda = \sum_{i=1}^{\ell}
k_i\,\alpha_i$ and $$\{\alpha_i\},\,i=1,\ldots,\ell,\eqno (2.12)$$
is the set of simple positive roots. 

Let $\hbox{min}\,m(\lambda)$ be the minimal value of $m_i(\lambda)$ for
$i=1,\ldots,m_{\ell(\lambda)}$. If $\ell(\lambda)> 1$, note that
$$\hbox{min}\,m(\lambda)<m_{\ell(\lambda)}(\lambda)\eqno(2.13)$$ by (2.10).

Clearly $H_{\lambda,i}(\lambda)$ is the highest
weight space of $H_{\lambda,i}$ and hence we also note, by (2.5), that
$$(S(\g)(\lambda))^{\n} = \sum_{i=1}^{\ell(\lambda)} J\otimes H_{\lambda,i}
(\lambda).\eqno
(2.14)$$\vskip.5pc Of course the left side of (2.14) is graded. It follows
immediately from (2.14) that if $S^k(\g)(\lambda)^{\n} \neq 0$, then
$$k\geq \hbox{min}\, m(\lambda).\eqno (2.15)$$ 

\vskip 1.5pc  {\bf 2.3.} The
universal enveloping algebra of a Lie algebra
$\ss$ is denoted by $U(\ss)$. Since we will be dealing with multiplication
in both
$U(\g)$ and $S(\g)$, when $x\in \g$, we will on
occasion to avoid confusion write $\tilde {x}$ for $x$ when $x$ is to be
regarded as an element of $U(\g)$ and not $S(\g)$. Of course $U(\g)$ is a
$G$-module by extension of the adjoint representation. By PBW one has a
$G$-module isomorphism $$\tau:S(\g)\to U(\g)\,\,\hbox{where for any $k\in
{\Bbb N}$ and $x\in \g$ one has $\tau(x^k) = \tilde {x}^k$}.\eqno (2.16)$$
\vskip .5pc 

Since
$\tau$ is a $\g$-module isomorphism the restriction of (2.16) to
$S(\g)^{\n}$ and to $S(\n)^{\n}$ readily yields linear isomorphisms
$$\tau:S(\g)^{\n}\to U(\g)^{\n}\eqno (2.17)$$ and $$\tau:S(\n)^{\n}\to
\hbox{\rm cent}\,U(\n).\eqno (2.18)$$ 

\vskip .5pc If
$\ss\s
\g$ is a Lie subalgebra and
$k\in \Bbb N$ the image, under
$\tau$, of $S^k(\ss)$, will be denoted by $U^{(k)}(\ss)$, and one readily
notes the direct sum $$ U(\ss) = \sum_{k=0}^{\infty}U^{(k)}(\ss).\eqno
(2.19)$$ If
$u\in U(\ss)$ then, by abuse of terminology, we will say that $u$ is
homogeneous of degree $k$ if $u\in U^{(k)}(\ss)$. We are particularly
interested in the case where $\ss = \n$ and $\n_-$. 

Since $\hh$ normalizes $\n$, both sides of (2.17) and (2.18) are bigraded
by degree and ($\hh$) weight and clearly (2.17) and (2.18) preserve the
bigrading. The proof of the following theorem uses results in [J]
where $\hbox{\rm cent}\,U(\n)$ is denoted by $Z(\n)$. 

\vs {\bf Theorem 2.1.}
(T. Joseph) {\it Let $\lambda\in \hbox{\rm Dom}(L_o)$ and assume that
$S(\n)^{\n}(\lambda)
\neq 0 $. Then  $S(\n)^{\n}(\lambda)$ is homogeneous of
 degree $\hbox{\rm min}\,m(\lambda)$.}

\vs {\bf Proof.} By (2.18) one has $\hbox{\rm cent}\,U(\n)(\lambda)\neq 0$. By
(iii), p. 260, in the Theorem of \S 4.12 of [J] and (iii), p.
261, in the lemma of \S 4.13 of [J] one has that
$\hbox{\rm cent}\,U(\n)(\lambda)$ is homogeneous of degree $j$ for some $j\in
\Bbb N$. But then
$S(\n)^{\n}(\lambda)$ is homogeneous of degree $j$ by (2.18). Consequently
$$j\geq \hbox{\rm min}\,m(\lambda)\eqno (2.20)$$ by (2.15). Assume $j>
 \hbox{min}\,m(\lambda)$. Then there exists $i \in \{1,\ldots, \ell(\lambda)\}$
such that
$$m_i(\lambda) < j.\eqno (2.21)$$ But $H_{\lambda,i}(\lambda)\s
(S^{m_i(\lambda)}(\g))^{\n}(\lambda)$. It follows then that
$(U^{(m_i(\lambda))}(\g))(\lambda)\neq 0$ by (2.17). But this and (2.21)
contradict  (iii), p.~261, in the lemma of \S 4.13 of [J].   Thus $j =
\hbox{min}\,m(\lambda)$.\hfill QED

\vskip 1pc {\bf 2.4.} Let $e =
\sum_{i=1}^{\ell}e_{\alpha_i}$ so that  $e$ is a principal nilpotent of
$\g$. Let $h\in\hh$ be fixed so that $\alpha_i(h) = 2$ for all simple
roots $\alpha_i$ so that $$[h,e]=2e.\eqno (2.22)$$ Then as one knows there
exists
$c_i\in \Bbb C^{\times},\,i=1,\ldots,\ell$, such that if
$e_-=\sum_{i=1}^{\ell}c_i\,e_{-\alpha_i}$, then $\{h,e,e_-\}$ is an
$Sl(2)$-triple and spans a principal TDS $\a$ in $\g$. Let $\xi\in
\hbox{\rm Dom}(L)$. Let $\xi^*\in \hbox{\rm Dom}(L)$ be the highest weight of
the $G$-module $V_{\xi}^*$ dual to
$V_{\xi}$. We retain this $*$ notation throughout. Then as $G$-modules 
$$V_{\xi}\otimes V_{\xi}^*\cong \hbox{\rm End}\,V_{\xi}.$$ Let $\nu =
\xi + \xi^*$ so that $V_{\nu}$ identifies with the Cartan product
of $V_{\xi}$ and $V_{\xi}^*$. Since the weights of
$V_{\xi}^*$ are the negatives of the weights of $V_{\xi}$ it
follows that $\nu\in \hbox{\rm Dom}(L_o)$. Furthermore since $\hbox{\rm
tr}\,AB$ defines a nonsingular invariant symmetric bilinear form on
$\hbox{\rm End}\,V_{\xi}$ it follows immediately that the corresponding
bilinear form on
$V_{\xi}\otimes V_{\xi}^*$ restricts to a nonsingular
$G$-invariant bilinear form
$(u,v)$ on $V_{\nu}$. Thus the highest and lowest weights in $V_{\nu}$
are respectively $\nu$ and $-\nu$. Consider the action of the
principal TDS $\a$ on $V_{\nu}$. Clearly, by the dominance of $\nu$ and
the regularity of $h$, the maximal (resp. minimal) eigenvalue of
$\pi_{\nu}(h)$ on
$V_{\nu}$ is
$\nu(h)$ (resp.
$-\nu(h)$ and these eigenvalues have multiplicity 1. Thus if $0 \neq
v_{\nu}$ and $0\neq v_{-\nu}$ are respectively highest and lowest 
 $\g$-weight vectors there
exists an irreducible $\a$-module $M$ in $V_{\nu}$ having $v_{\nu}$ and $v_{-\nu}$ respectively
as highest ($\Bbb C\, h$) weight vectors. Hence for some $z\in \Bbb C^{\times}$ one has
$$(\pi_{\nu}(e_-))^{\nu(h)}v_{\nu} = z\,v_{-\nu}.\eqno (2.23)$$ Let
$\lambda = 2\nu$ so that certainly $\lambda\in \hbox{\rm Dom}\,L_o$.  Recalling
(2.11) and the 2 in (2.22) one has $$\nu(h) =
m_{\ell(\lambda)}(\lambda).\eqno (2.24)$$ But clearly 
$(v_{-\nu},v_{\nu})\neq 0$ (because of
multiplicity 1). Thus one has $$(\pi_{\nu}({\widetilde
e_-}^{\,\,\,\,\,\,m_{\ell(\lambda)}(\lambda)})v_{\nu},v_{\nu})\neq 0.\eqno
(2.25)$$ For the tilde notation see (2.16).

In the next section we will give a simple condition guaranteeing that
$\ell(\lambda)>1$.

\vskip 1.5pc {\bf 2.5.}  We noted the following result,
Proposition 2.2, a long time ago. However it is likely that the result was
well known even then but we are unable to find published references to it
so, for completeness, we will give a proof here. The applications, Theorems
2.4 and 2.5, of Proposition 2.2, are recent with me.  Let
$\beta,
\gamma\in \hbox{\rm Dom}(L)$. Then, as one knows, the Cartan product $V_{\beta
+
\gamma}$ occurs with multiplicity one in the tensor product
$V_{\beta}\otimes V_{\gamma}$ so that there exists a unique $G$-invariant
projection 
$$\Gamma: V_{\beta}\otimes V_{\gamma}\to V_{\beta + \gamma}.\eqno
(2.26)$$

\vskip .5pc {\bf Proposition 2.2.} {\it Let $0\neq u\in
V_{\beta}$ and $0\neq w\in V_{\gamma}$. Then $$\Gamma(u\otimes w)\neq 0.
\eqno (2.27)$$}
 \indent {\bf Proof.} Let $0\neq v_{\beta}$ (resp.
$0\neq v_{\gamma})$ be a highest weight vector in $V_{\beta}$ (resp.
$V_{\gamma})$ so that $v_{\beta}\otimes v_{\gamma}$ is a highest weight
vector in the Cartan product $V_{\beta + \gamma}$. Consequently taking
into account the action of $G$,
$$\{g\cdot v_{\beta}\otimes g\cdot v_{\gamma}\mid g\in
G\}\,\,\,\,\hbox{spans}\,\,\,  V_{\beta + \gamma}.\eqno(2.28)$$\vskip .2pc

Let
$K$ be a maximal compact subgroup of $G$ (so that $G$ is the
complexification of $K$) and let $\{y,z\}_{\beta}$
 (resp. $\{y,z\}_{\gamma}$) be a $K$-invariant Hilbert space structure on
$V_{\beta}$ (resp. $V_{\gamma}$). These induce a natural $K$-invariant
Hilbert space structure $\{y,z\}_{\beta,\gamma}$ on $V_{\beta}\otimes
V_{\gamma}$. Furthermore it is immediate that $\Gamma$ is a Hermitian
projection with respect to the latter inner product. Thus to prove the
proposition it suffices to show that there exists $g_o\in G$
such that $$\{g_o\cdot v_{\beta}\otimes g_o\cdot v_{\gamma}, u\otimes
w\}_{\beta,\gamma}\neq 0.\eqno (2.29)$$ But now since $V_{\beta}$ (resp.
$V_{\gamma}$) is the span of $\{g\cdot v_{\beta}\mid g\in G\}$ (resp.
$\{g\cdot v_{\gamma}\mid g\in G\}$) it follows immediately that the
function $F_{\beta}$ (resp. $F_{\gamma}$) on $G$ given by $F_{\beta}(g) =
\{g\cdot v_{\beta},u\}_{\beta}$ (resp. $F_{\gamma}(g) =
\{g\cdot v_{\gamma},w\}_{\gamma}$ ) is nonvanishing and analytic. Thus
there exists $g_o\in G$ such that $F_{\beta}(g_o)F_{\gamma}(g_o)\neq 0$.
But the left side of (2.29) equals $F_{\beta}(g_o)F_{\gamma}(g_o)$. This
proves (2.29).

\hfill QED

\vs As a corollary one has \vs {\bf Proposition 2.3.}
{\it Let the notation be as in Proposition 2.2. Assume that $\ss$ is any
subspace of $V_{\beta}$. Then the map $$\ss\to V_{\beta +
\gamma},\,\,\,x\mapsto \,\Gamma(x\otimes w)\eqno (2.30)$$ is
linear and injective. Furthermore if $\ss\s V_{\beta}(\mu)$ for some
$\mu\in L$ and $w\in V_{\gamma}(\delta)$, for some $\delta\in L$, then the
image of (2.30) is contained in $V_{\beta +\gamma}(\mu + \delta)$. }

 \vs
{\bf Proof.} Obviously (2.30) is linear. The injectivity is
immediate from Proposition 2.2. The second conclusion follows from the
fact that 
$\Gamma$ is, among other things, an $\hh$-map.\hfill QED

 \vs We now have the
following information about a $0$-weight space.  

\vs {\bf Theorem 2.4.}
{\it Let $\beta\in \hbox{\rm Dom}(L)$ and let $\beta^*$ be the highest weight
of the contragredient module to $V_{\beta}$. Let
$d$ be the maximal value of all the multiplicities of weights in
$V_{\beta}$. Then
$$\hbox{\rm dim}\, V_{\beta + \beta^*}(0) \geq d. \eqno (2.31)$$}

\vs {\bf
Proof.} We retain the notation of Propositions 2.2 and 2.3. Choose
$\gamma= \beta^*$ and let $\mu$ be any weight of $V_{\beta}$. Choose 
$\delta = -\mu$. Then the image of (2.30) is contained in $V_{\beta +
\beta^*}(0)$ by Proposition 2.3. But the dimension of the image equals
$\hbox{\rm dim}\,V_{\beta}(\mu)$ by the injectivity of (2.30). Since $\mu$ is
arbitrary this of course implies (2.31).\hfill QED

\vs A similar argument leads
to the following monotonicity result of weight
multiplicities.

 \vs {\bf Theorem 2.5.} {\it Let $\beta\in \hbox{\rm Dom}(L)$ and
$\gamma \in \hbox{\rm Dom}(L_o)$. Then for any weight $\mu$ of $V_{\beta}$ one
has 
$$\hbox{\rm dim}\,V_{\beta + \gamma}(\mu)\geq \hbox{\rm
dim}\,V_{\beta}(\mu).\eqno (2.32)$$} 

{\bf Proof.} Again we use the notation and
result in Propositions 2.3.
 Since
$\gamma\in \hbox{\rm Dom}(L_o)$ one knows $V_{\gamma}(0)\neq 0$. Let $\delta
=0$ in Proposition 2.3 and let $\ss = V_{\beta}(\mu)$ so that the image of
(2.30) is contained in $V_{\beta + \gamma}(\mu)$. But then (2.32) follows
from the injectivity of (2.30).\hfill QED

\vs We return to the notation of \S
2.4. Recall that $\xi\in \hbox{\rm Dom}(L),\,\, \nu = \xi + \xi^*$ and $\lambda
= 2\,\nu$ so that $\nu,\,\lambda$ are in $\hbox{\rm Dom}(L_o)$. Let
$d$ be the maximal value of all weight multiplicities of $V_{\xi}$.

\vs {\bf
Theorem 2.6.} {\it Assume $$d>1.\eqno (2.33)$$ Then
$\hbox{\rm dim}\,V_{\lambda}(0)>1$ so that $$\hbox{\rm min} \,m(\lambda)<
m_{\ell(\lambda}(\lambda). \eqno (2.34)$$ Furthermore $\xi$ can be chosen so
that (2.33) is satisfied if and only if there exists a simple component of
$\g$ which is not of type $A_1$. }

\vs {\bf Proof.} One has $\hbox{\rm dim}\,V_{\nu}(0)
>1$ by Theorem 2.4. But then
$\hbox{\rm dim}\,V_{\lambda}(0) >1 $ by Theorem 2.5. The statement (2.34) is
then just (2.13). If all the simple components of $\g$ are of type $A_1$,  then
clearly
$d=1$ for any $\xi\in \hbox{\rm Dom}(L)$. However, if not, then $d>1$ for the
adjoint representation of a component, not of type of $A_1$, when extended
trivially to the other components. \hfill QED

\vskip 1.5pc {\bf 2.6.} In [LW]
the authors introduce a very neat idea for constructing certain elements in
$S(\n)^{\n}$ (or equivalently in $\hbox{\rm cent}\,U(\n)$) using representation
theory. The statement of this idea, Lemma 3.7 in [LW], however, is not
correct as it stands. We give a counterexample in this section. Nevertheless
the statement of this lemma in [LW] can be modified so as to establish that
this very interesting technique in [LW] does indeed yield elements in
$S(\n)^{\n}$. We do this in the section that follows this one. 

The bilinear form ${\cal B}$ on $S(\g)$ clearly defines a nonsingular
pairing of $S(\n_-)$ and $S(\n)$ with $S^i(\n_-)$ orthogonal to $S^j(\n)$
when $i\neq j$ and $S^i(\n)\cong  S^i(\n_-)^*$.  Let the following notation be
as in \S 2.4 and let $f_{\nu}\in  S(\n) $ be defined so that if $\Xi\in
S(\n_-)$, then $$(f_{\nu}, \Xi) = (\pi_{\nu}(\tau(\Xi))\,
v_{\nu},v_{\nu})\eqno(2.35)$$ where $\tau$ is defined as in (2.16). One notes
that (5.35) vanishes if $\Xi\in S^j(\n_-)$ for $j> \hbox{\rm dim}\,V_{\nu}$ so
that 
$f_{\nu}\in  S(\n) $ is well defined. Lemma 3.7 in [LW] asserts that
$$f_{\nu}\in S(\n)^{\n}(\lambda). \eqno (2.36)$$ But (2.35) is not $0$ if
$\Xi = (e_-)^{m_{\ell(\lambda)}(\lambda)}$ by (2.25). On the other hand
$S(\n)^{\n}(\lambda)\s S^{\hbox{\sevenrm min}\,m(\lambda)}(\g)$ by Theorem 2.1
and
$(e_-)^{m_{\ell(\lambda)}(\lambda)}\in S^{m_{\ell(\lambda)}(\lambda)}(\g)$.
But one has (2.34) if
$\xi$ and $\g$ are chosen as in the last statement in Theorem 2.6. Such a
choice leads to a contradiction since (2.34) implies
$S^{\hbox{\sevenrm min}\,m(\lambda)}(\g)$ is orthogonal to
$S^{m_{\ell(\lambda)}(\lambda)}(\g)$.

\vskip 1pc {\bf 2.7.} Let $U_n(\g),
\,n\in \Bbb N$, be the standard filtration of $U(\g)$. A nonzero linear
functional $f$ on $U(\g)$ will be said to have codegree $k$ if $k\in \Bbb N$ is
maximal such that $f$ vanishes on $U_{k-1}(\g)$ (putting $U_{-1}(\g) = 0$).
Assume $f$ has codegree $k$ and $k\geq 1$. Note that if $x_i\in\g,\,
i=1,\ldots, k$, then for any permutation $\sigma$ of $\{1,\ldots,k\}$ one has
$$f(\widetilde {x}_1\cdots \widetilde {x}_k) = 
f(\widetilde{x}_{\sigma(1)}\cdots
\widetilde{x}_{\sigma(k)}) \eqno (2.37)$$ using the notation of (2.16). \vskip
.5pc

 Now let
$f_{(k)}$ be the linear functional on $S^k(\g)$ defined by the restriction
$f\circ
\tau\mid S^k(\g)$. Using ${\cal B}$ on $S(\g)$ we may regard $f_{(k)}$ as an
element in
$S^k(\g)$ so that by (2.37) one has $$f(\widetilde {x}_1\cdots \widetilde
{x}_k) = (f_{(k)},x_1\cdots x_k)\eqno (2.38)$$ for any $x_i\in\g,\,
i=1,\ldots, k$. One also notes that if $k = i +j$, where $i,j\in \Bbb N$, and
$v\in S^i(\g),\,w\in S^j(\g)$, then $$f(\tau(v)\tau(w)) = (f_{(k)}, v\,w)\eqno
(2.39)$$ since, clearly $\tau(v)\tau(w) - \tau(v\,w) \in U_{k-1}(\g)$. 

Taking a clue from [LW] we will choose $f$ so that it arises from a matrix
entry of a $U(\g)$-module. In fact assume $M$ is a $U(\g)$-module (not
necessarily finite dimensional) with respect to a representation $\pi:U(\g)\to
 \hbox{\rm End}\,M$. We recall that $\pi$ has an infinitesimal character
$\chi:\hbox{\rm cent}\,U(\g)\to \Bbb C$ if for any $z\in \hbox{\rm cent},U(\g)$
one has
$\pi(z) =
\chi(z)\,\hbox{\rm Id}_M$. Then one has

 \vs {\bf Theorem 2.7.} {\it Assume
$\pi$ is a representation of $U(\g)$ on a vector space $M$ and that $\pi$ has
an infinitesimal character $\chi$. Assume also if $s\in M$ and $s'\in M^*$ are
such that $f\in U(\g)^*$,  defined by $f(u) = s'(\pi(u)s)$, for any
$u\in U(\g)$, is nonvanishing with codegree $k\geq 1$. Then
$f_{(k)}$ is harmonic. That is $f_{(k)}\in H^k$.}

\vs {\bf Proof.} Let $i\in
\Bbb N$ where $i\geq 1$. We must show that if $r\in S^i(\g)^{\g}$, then
$$\partial_r f_{(k)}= 0.\eqno (2.40)$$ Obviously one has (2.40) if $i>k$ so
that we can assume that $i
\leq k$. Let $j = k-i$, so that $j<k$, and let $w\in S^j(\g)$ be arbitrary so
that it suffices to show that $$(\partial_r f_{(k)},w)= 0.\eqno (2.41)$$ But by
(2.1) the left side of (2.41) equals $(f_{(k)}, r\,w)$. But clearly
$\tau(r)\in \hbox{\rm cent}\,U(\g)$. Let $c= \chi(\tau(r))$. Then by (2.39)
$$\eqalign{(f_{(k)}, r\,w) &= f(\tau(r)\tau(w))\cr &=
c\,f(\tau(w))\cr}.$$ But $\tau(w)\in U_{k-1}(\g)$ since $j<k$. Hence 
$f(\tau(w)) = 0$. Thus $(f_{(k)}, r\,w)= 0$. This establishes (2.40).
\hfill QED\vs 

The main result, to follow, will use notations of \S 2.4 but with
basically fewer restrictive conditions. We eliminate $\xi$ and $0\neq \nu\in
\hbox{\rm Dom}(L)$ is now arbitrary. Put $\lambda = \nu + \nu^*$. Also
$V_{\nu^*}$ will be identified with the dual space $V_{\nu}^*$ to $V_{\nu} $.
Let $0\neq v_{\nu}$ (resp. $0\neq v_{\nu^*}$) be a highest weight vector of
$V_{\nu}$ (resp. $V_{\nu^*}$). Let
$f\in U(\g)^*$ be defined by $$f(u)= v_{\nu^*}(\pi_{\nu}(u)v_{\nu}) \eqno
(2.42)$$ for any $u\in U(\g)$. Expressed another way regard $v_{\nu}\otimes
v_{\nu^*}$ as a rank 1 linear operator on $V_{\nu}$ where for any $v\in V_{\nu}$
one has $v_{\nu}\otimes v_{\nu^*}(v) = v_{\nu^*}(v)\,v_{\nu}$. Then $$f(u) =
tr\, \pi_{\nu}(u)\,\,v_{\nu}\otimes v_{\nu^*}.\eqno (2.43)$$ Now let $k$ be the
codegree of $f$. But $v_{\nu^*}(v_{\nu}) = 0$ since $0\neq \nu$. Thus $k\geq
1$. It is then immediate from (2.19) and (2.43) that $$f_{(k)}\in
S(\g)^{\n}(\lambda).\eqno (2.44)$$ The main point is to show that we may
replace $\g$ by $\n$ in (2.44). 

Let $\u = S^k(\g)\cap (\hh + \n)S(\g)$. Clearly one has a direct sum
$$S^k(\g) = S^k(\n_-) \oplus \u.\eqno (2.45)$$ Furthermore it is immediate that
$\u$ is the ${\cal B}$ orthocomplement of $S^k(\n)$ in $S^k(\g)$. Thus to prove
that $$f_{(k)}\in S(\n)^{\n}(\lambda)\eqno (2.46)$$ it suffices to show that
$f_{(k)}$ is
${\cal B}$ orthogonal to $\u$. Clearly any element in $\u$ is a sum of elements
of the form $y\,w$ where $w\in \hh + \n$ and $y\in S^{k-1}(\g)$. However given
such an element there exists $a\in \Bbb C$ such that $\pi_{\nu}(\tau(w))v_{\nu} =
a\,v_{\nu}$. Also $f(\tau(y))= 0$ since $f$ vanishes on $U_{k-1}(\g)$. But then
by (2.39)
$$\eqalign{(f_{(k)}, y\,w)&= f(\tau(y)\tau(w))\cr &=
v_{\nu^*}(\pi_{\nu}(\tau(y)\tau(w))v_{\nu}) \cr &= a\,
v_{\nu^*}(\pi_{\nu}(\tau(y))v_{\nu})\cr &= a\,f(\tau(y))\cr &= 0.\cr}$$ This
proves (2.46). That is, we have proved 

\vs {\bf Theorem 2.8.} {\it Let
$0\neq \nu\in \hbox{\rm Dom}(L)$ and let $\lambda= \nu +\nu^*$. Let $f\in U(\g)^*$
be defined by (2.42) or equivalently (2.43). Let $k$ be the codegree of $f$. Then
$k\geq 1$ and
$$f_{(k)}\in S(\n)^{\n}(\lambda).\eqno (2.47)$$}

\vskip .2pc

\centerline{\bf References}\vskip 1pc
\item {[LW]} Ronald Lipsman and Joseph Wolf, {\it Canonical
semi-invariants and the Plancherel formula for parabolic groups},
Trans. Amer. Math. Soc., {\bf 269}(1982), 111--131.
\item {[J]} Anthony Joseph, {\it A preparation theorem for the
prime spectrum of a semisimple Lie algebra}, Journ. of Alg., {\bf
48}, No.2, (1977), 241--289.
\item {[K]} Bertram Kostant, {\it The cascade of orthogonal roots
and the coadjoint action of a maximal unipotent subgroup of a
semisimple Lie group}, to appear.

\vs

\noindent Bertram Kostant (Emeritus)

\noindent Department of Mathematics

\noindent Cambridge, MA 02139

\noindent email: kostant@math.mit.edu

\end